\documentclass{amsart}

\usepackage{amssymb}

\usepackage{graphicx}

\usepackage[cmtip,all]{xy}

\copyrightinfo{2006}{American Mathematical Society}

\newtheorem{theorem}{Theorem}[section]
\newtheorem{lemma}[theorem]{Lemma}
\newtheorem{cor}[theorem]{Corollary}
\newtheorem{prop}[theorem]{Proposition}

\theoremstyle{definition}

\newcommand{\Nbd}{\operatorname{Nbd}}
\newcommand{\cl}{\operatorname{cl}}

\theoremstyle{remark}

\numberwithin{equation}{section}

\begin{document}

\title[Homeomorphisms of the 3-sphere]{Homeomorphisms of the 3-sphere that preserve a Heegaard
splitting of genus two}

\author{Sangbum Cho}
\address{Department of Mathematics, University of Oklahoma, Norman, Oklahoma 73019}
\email{scho@math.ou.edu}


\date{April 15, 2007}

\begin{abstract}
Let ${\mathcal H_2}$ be the group of isotopy classes of
orientation-preserving homeomorphisms of $\Bbb S^3$ that preserve
a Heegaard splitting of genus two. In this paper, we construct a
tree in the barycentric subdivision of the disk complex of a
handlebody of the splitting to obtain a finite presentation of
${\mathcal H_2}$.
\end{abstract}

\maketitle

\section{Introduction}

Let $\mathcal H_g$ be the group of isotopy classes of
orientation-preserving homeomorphisms of $\Bbb S^3$ that preserve a
Heegaard splitting of genus $g$, for $g \geq 2$. It was shown by
Goeritz \cite{Go} in 1933 that $\mathcal H_2$ is finitely generated.
Scharlemann \cite{Sc} gave a modern proof of Goeritz's result, and
Akbas \cite{Ak} refined this argument to give a finite presentation
of $\mathcal H_2$. In arbitrary genus, first Powell \cite{Po} and
then Hirose \cite{Hi} claimed to have found a finite generating set
for the group $\mathcal H_g$, though serious gaps in both arguments
were found by Scharlemann. Establishing the existence of such
generating sets appears to be an open problem.

In this paper, we recover Akbas's presentation of $\mathcal H_2$ by
a new argument. First, we define the complex $P(V)$ of primitive
disks, which is a subcomplex of the disk complex of a handlebody $V$
in a Heegaard splitting of genus two. Then we find a suitable tree
$T$, on which $\mathcal H_2$ acts, in the barycentric subdivision of
$P(V)$ to get a finite presentation of ${\mathcal H_2}$. In the last
section, we will see that the tree $T$ can be identified with the
tree used in Akbas's argument \cite{Ak}.

Throughout the paper, $(V, W; \Sigma)$ will denote a Heegaard
splitting of genus two of $\Bbb S^3$. That is, $\Bbb S^3 = V \cup W
$ and $ V \cap W =\partial V =\partial W = \Sigma $, where $V$ and
$W$ are handlebodies of genus two. For essential disks $D$ and $E$
in a handlebody, the intersection $D \cap E$ is always assumed to be
transverse and minimal up to isotopy. In particular, if $D$
intersects $E$ (indicated as $D \cap E \neq \emptyset$), then $D
\cap E$ is a collection of pairwise disjoint arcs that are properly
embedded in both $D$ and $E$.

Finally, $\Nbd(X)$ will denote a regular neighborhood of $X$ and
$\cl(X)$ the closure of $X$ for a subspace $X$ of a polyhedral
space, where the ambient space will always be clear from the
context.

\smallskip

\begin{center}
{\sc Acknowledgement}
\end{center}

\smallskip

The author would like to thank his advisor D. McCullough for his
consistent encouragement and sharing his enlightening ideas on the
foundations of this project. The author also would like to thank M.
Scharlemann for his valuable suggestions and corrections. Finally,
the author is grateful to the referee for carefully reading the
paper and suggesting many improvements.

\begin{figure}
\centering
\includegraphics{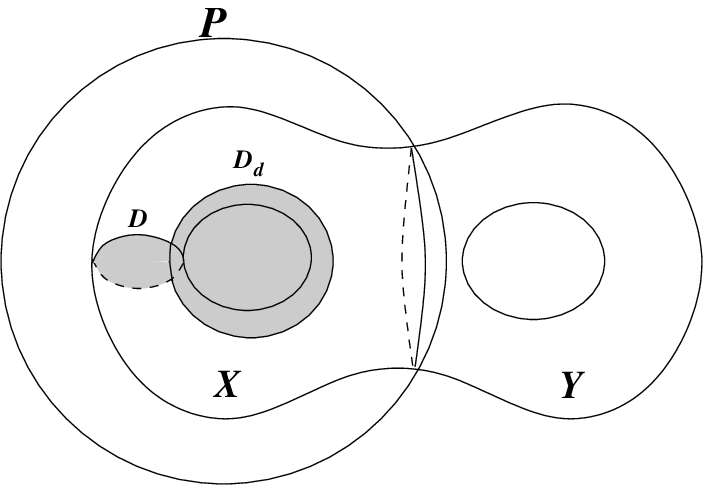}
\caption{}
\end{figure}

\section{Primitive disks in a handlebody}

We call an essential disk $D$ in $V$ {\it primitive} if there exists
an essential disk $D_d $ in $W$ such that $\partial D$ and $\partial
D_d $ have a single transverse intersection in $\Sigma$. Such a
$D_d$ is called a {\it dual disk} of $D$. Notice that any primitive
disk is nonseparating. We call a pair of disjoint, nonisotopic
primitive disks in $V$ a {\it reducing pair} of $V$. Similarly, a
triple of pairwise disjoint, nonisotopic primitive disks is a {\it
reducing triple}.

A 2-sphere $P$ in $\Bbb S^3$ is called a {\it reducing sphere} for
$(V,W;\Sigma)$ if $P$ intersects $\Sigma$ transversally in a
single essential circle and so intersects each handlebody in a
single essential disk. It is clear that $V \cap P$ and $ W \cap P$
are essential separating disks in $V$ and $W$ respectively.

\begin{lemma}
A disk $D$ in $V$ is primitive if and only if $\cl(V-\Nbd(D))$ is
an unknotted solid torus, i.e. $W \cup \Nbd(D)$  is a solid torus.
\end{lemma}

\begin{proof}
 Let $D$ be a primitive disk with a dual disk $D_d$.
 Then the boundary $P$ of a regular neighborhood of $D \cup D_d$ is
a reducing sphere that is disjoint from $D$ (see Fig. 1).
 The sphere $P$ splits $V$ into two solid tori, $X$ and $Y$, and we may assume
that $D$ is a meridian disk for $X$.
 Cutting $\Bbb S^3$ along $P$, we get two 3-balls, $B_1$ and $B_2$,
where $X$ and $Y$ are contained in $B_1$ and $B_2$ respectively.
 Since the handlebody $W$ is the boundary connected sum of
$\cl(B_1-X)$ and $\cl(B_2-Y)$ along the disk $W \cap P$, we have
$$\Bbb Z * \Bbb Z = \pi_1 (W) = \pi_1(B_1 - X) * \pi_1(B_2 -
Y) =  \pi_1 (\Bbb S^{3} -X)* \pi_1 (\Bbb S^{3} -Y).$$
 Thus $ \pi_1 (\Bbb S^{3} -X)= \pi_1 (\Bbb S^{3} -Y)= \Bbb Z$,
and consequently $X$ and $Y$ are unknotted.
 Since $D$ is a meridian of $X$, so $\cl(V-\Nbd(D))$ is ambient
isotopic to $Y$, and is thus unknotted.
 The converse is a special case of Theorem 1 in \cite{Go}.
\end{proof}

 Let $E$ and $D$ be nonseparating disks in $V$ such that
$E \cap D \neq \emptyset$, and let $C \subset D$ be a disk cut off
from $D$ by an outermost arc $\alpha$ of $D \cap E$ in $D$ such
that $C\cap E= \alpha$.
 The arc $\alpha$ cuts $E$ into two disk components, say $G$ and $H$.
 The two disks $E_1 = G \cup C$ and $E_2 = H \cup C$ are called the
{\it disks obtained from surgery on $E$ along $C$}.
 Since $E$ and $D$ are assumed to intersect minimally, $E_1$
and $E_2$ are isotopic to neither $E$ nor $D$, and moreover have
fewer arcs of intersection with $D$ than $E$ had.

 Notice that $E_1$ and $E_2$ are not isotopic to each other,
otherwise $E$ would be a separating disk.
 Finally, observe that both $E_1$ and $E_2$ are isotopic
to a meridian disk of the solid torus $\cl(V - \Nbd(E))$, otherwise
one of them would be isotopic to $E$.
 Thus $E_1$ and $E_2$ are all nonseparating disks in
$V$.

\begin{lemma}
For a reducing pair \{$E, E'$\} of  $V$, there exists, up to
isotopy, a unique reducing pair $\{E_d, E'_d \}$ of  $W$ for which
$E_d$ and $E'_d$ are dual disks of $E$ and $E'$, and are disjoint
from $E'$ and $E$ respectively.
\end{lemma}

\begin{proof}
 Existence is also a special case of Theorem 1 in \cite{Go}.
 For uniqueness, assume that \{$F_d , F'_d$\} is
another reducing pair of $W$ such that $F_d$ and $F'_d$ are dual
disks of $E$ and $E'$, and are disjoint from $E'$ and $E$
respectively.

 Suppose that $F_d$ intersects $E'_d$.
 Then there exists a disk $C\subset F_d$ cut off from $F_d$ by an
outermost arc $\alpha$ of $F_d \cap E'_d$ in $F_d$ such that
$C\cap E'_d= \alpha$.
 Denote the arc $C\cap\Sigma$ by $\delta$.
 The disks obtained from surgery on $E'_d$ along $C$ are isotopic to a meridian
disk of the solid torus $\cl(W - \Nbd(E'_d))$.
 Hence $\delta$ must intersect $\partial E$ which is a longitudinal
circle of the solid torus.

 Let $\tilde \Sigma$ be the 2-holed torus obtained by cutting $\Sigma$
along $E'_d$.
 Denote by $l_\pm$ the boundary circles of $\tilde \Sigma$ that came
from $\partial E'_d$.
 Then $\delta$ is an essential arc in $\tilde\Sigma$ with endpoints
in a single boundary circle $l_+$ or $l_-$, say $l_+$.
 Moreover, since $\partial F_d \cap l_+$ and $\partial F_d \cap
l_-$ contain the same number of points, we also have an essential
arc $\delta'$ of $\partial F_d \cap\tilde \Sigma $ with endpoints
in the boundary circle $l_-$. The arc $\delta'$ also intersects
$\partial E$, otherwise $\delta'$ would be inessential.
 Thus $\partial F_d$ meets $\partial E$ in at least two points and
this contradicts that $F_d$ is a dual disk of $E$.
 Therefore, we conclude that $F_d$ is disjoint from $E'_d$.

 Now, both $E_d$ and $F_d$ are meridian disks of the
solid torus $\cl(W - \Nbd(E'_d))$.
 Since both $E_d$ and $F_d$ are disjoint from $E'\cup E'_d$, they are
isotopic to each other in the solid torus and hence in $W$.
Similarly, $F'_d$ is isotopic to $E'_d$ in $W$.
\end{proof}

\begin{theorem}
Let $E$ and $D$ be primitive disks in $V$ such that $E \cap D \neq
\emptyset$ and let $C \subset D$ be a disk cut off from $D$ by an
outermost arc $\alpha$ of $D \cap E$ in $D$ such that $C\cap E=
\alpha$. Then both disks obtained from surgery on $E$ along $C$
are primitive.
\end{theorem}

\begin{proof}
 First, we consider the case when there exists a primitive disk
$E'$ such that $\{E, E'\}$ is a reducing pair and $E'$ is disjoint
from $C$.
 (The existence of such an $E'$ will be established in the final
 paragraph.)
 Let $E_d$ and $E_d'$ be the dual disks of $E$ and $E'$ respectively
given by Lemma 2.2.
 By isotopy of $D$, we may assume that $\partial D$ intersects
$\partial E_d$ and $\partial E'_d $ minimally in $\Sigma$.
 Denote the arc $C \cap \Sigma$ by $\delta$.
 It suffices to show that $\delta$ intersects $\partial E_d '$ in a
single point, since then the resulting disks from surgery are both
primitive with common dual disk $E_d '$.

 Let $\Sigma'$ be the 4-holed sphere obtained by cutting $ \Sigma$
along $\partial E \cup \partial E'$.
 Denote by $\partial E_\pm$ (resp. $\partial E'_\pm$) the boundary
circles of $\Sigma'$ that came from $\partial E$ (resp. $
\partial E'$).
 Then $ \delta $ is an essential arc in $\Sigma'$ with endpoints in a
single boundary circle $\partial E_+$ or $\partial E_-$, say
$\partial E_-$, and $\delta$ cuts off an annulus from $\Sigma'$.
 The boundary circle of the annulus that does not contain $\delta$
cannot be $\partial E_+$ otherwise one of the disks obtained from
surgery on $E$ along $C$ would be isotopic to $E$.
 Thus it must be $\partial E'_+$ or $\partial E'_-$, say $\partial E'_+$.
 Let $\gamma$ be a spanning arc of the annulus connecting $\partial
E_-$ and $\partial E'_+ $.
 Then $\delta$ can be regarded as the frontier of
$\Nbd (\partial E'_+ \cup \gamma)$ in $\Sigma'$ (see Fig. 2(a)).
 In $\Sigma'$, the boundary circle $\partial E_d$ (resp. $\partial E'_d $)
appears as an arc connecting $\partial E_+$ and $\partial E_-$
(resp. $\partial E'_+$ and $\partial E'_- $).
 We observe that $\delta$ intersects $\partial E_d '$ in at least
 one point.

\begin{figure}
\centering
\includegraphics{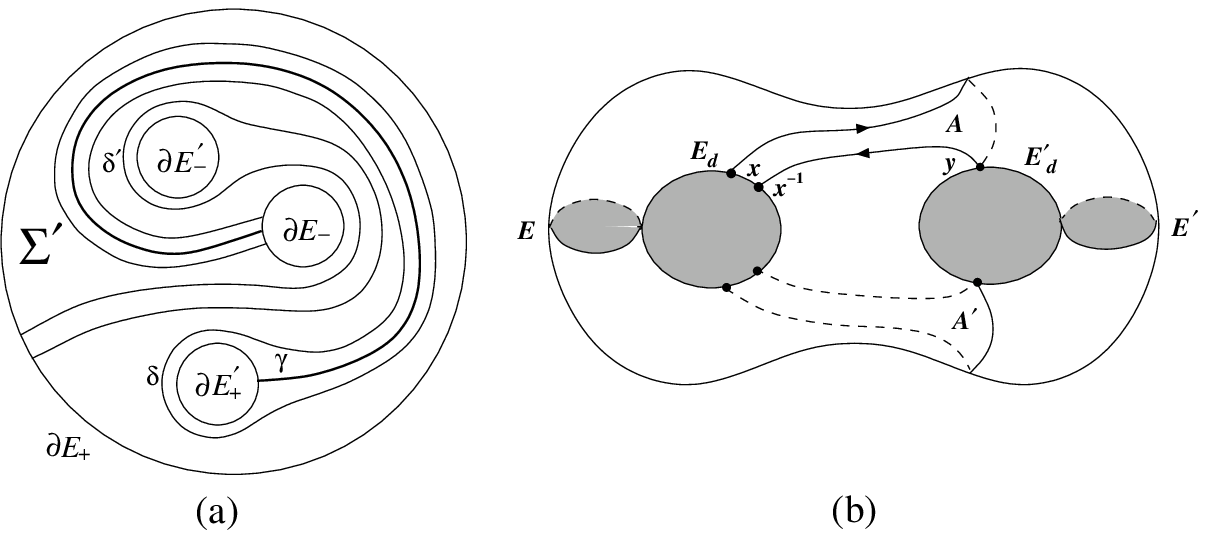}
\caption{}
\end{figure}

 Suppose, for contradiction, that $\delta$ intersects $\partial E'_d$ in more than one point.
 Then $\gamma$ intersects $\partial E_d$ (and $\partial E'_d$)
at least once. In particular, there exists an arc component of
$\gamma \cap (\Sigma' -\partial E_d)$, which connects $\partial
E'_+$ and $\partial E_d$.
 Consequently, $\delta$ cuts off an annulus $A$ from
$\Sigma' -\partial E_d$, having one boundary circle $\partial E'_+$
(see Fig. 2(b)).

 Now, $\partial E$ and $\partial E'$ represent the two generators
$x$ and $y$, respectively, of the free group $\pi_1 (W) = \langle
x, y\rangle$, and $\partial D$ represents a word $w$ in terms of
$x$ and $y$.
 Each such word $w$ can be read off from the intersections of
$\partial D$ with $\partial E_d$ and $\partial E'_d$.
 In particular, the arc $\delta \cap \partial A$ determines
a sub-word of a word $w$ of the form $xy^{\pm 1} x^{-1}$ or
$x^{-1}y^{\pm 1} x$, and hence each word $w$ contains both $x$ and
$x^{-1}$ (see Fig. 2(b)).

 We claim that each word $w$ is reduced, and therefore cyclically
reduced.
 Since $\partial D \cap \partial E_+$ and $\partial D \cap \partial E_-$
contain the same number of points, we also have an essential arc
$\delta'$ of $\partial D \cap \Sigma'$ whose endpoints lie in
$\partial E_+$. The arc $\delta'$ cuts off an annulus from
$\Sigma'$, and the boundary circle of the annulus that does not
contain $\delta'$ must be $\partial E'_-$ (see Fig. 2(a)). Since
$\delta'$ also intersects $\partial E'_d$ in more than one point,
the above argument holds for $\delta'$.
 In particular, $\delta'$ cuts off an annulus $A'$ from $\Sigma'
-\partial E_d$ having one boundary circle $\partial E'_-$.
 Observe that the annuli $A$ and $A'$ meet $\partial E_d$ on opposite
sides from each other, as in Fig. 2(b), since otherwise $\partial D$
would not have minimal intersection with $\partial E_d \cup
\partial E'_d$.

 Let $\Sigma''$ be the 4-holed sphere obtained by cutting $\Sigma$
along $\partial E_d \cup \partial E'_d$.
 Then $\partial A \cup \partial A'$ cuts $\Sigma''$ into two disks,
and each disk meets each boundary circle of $\Sigma''$ in a single
arc.  Consequently, we see there exists no arc component of
$\partial D$ in $\Sigma''$ that meets only one of $\partial E_d$
and $\partial E'_d$ in the same side.
 Thus we conclude that $w$ contains neither $x^{\pm1} x^{\mp1}$
nor $y^{\pm1} y^{\mp1} $.
 Since this is true for each word $w$, so each is cyclically reduced.

 It is well known that a cyclically reduced word $w$ in the free group
$\langle x, y\rangle$ of rank two cannot be a generator if $w$
contains $x$ and $x^{-1}$ simultaneously.
 Therefore, $\pi_1 (W\cup\Nbd(D)) =\langle x, y ~|~ w \rangle $ cannot
be the infinite cyclic group, and consequently $W \cup \Nbd(D)$ is
not a solid torus. This contradicts, by Lemma 2.1, that $D$ is
primitive in $V$.

It remains to show that such a primitive disk $E'$ does exist.
Choose a primitive disk $E'$ so that $\{E, E' \}$ is a reducing
pair. If $C$ is disjoint from $E'$, we are done. Thus suppose $C$
intersects $E'$. Then we have a disk $F \subset C$ cut off from $C$
by an outermost arc $\beta$ of $C\cap E'$ in $C$ such that $F \cap
E' = \beta$. By the above argument, the disks obtained from surgery
on $E'$ along $F$ are primitive. One of them, say $E''$, forms a
reducing pair $\{ E, E''\}$ with $E$. It is clear that $| C \cap E''
| < |C \cap E' | $ since at least $\beta$ no longer counts. Then
repeating the process for finding $E''$, we get the desired
primitive disk.
\end{proof}

\section{A sufficient condition for contractibility}

In this section, we introduce a sufficient condition for
contractibility of certain simplicial complexes. Let $K$ be a
simplicial complex. A vertex $w$ is said to be {\it adjacent} to a
vertex $v$ if equal to $v$  or if $w$ spans a 1-simplex in $K$ with
$v$. The {\it star} $st(v)$ of $v$ is the maximal subcomplex spanned
by all vertices adjacent to $v$.

A {\it multiset} is a pair $(A, m)$, typically abbreviated to $A$,
where $A$ is a set and $m:A \rightarrow \Bbb N$ is a function. In
other words, a multiset is a set with multiplicity. An {\it
adjacency pair} is a pair $(X, v)$, where $X$ is a finite multiset
whose elements are vertices of $K$ which are adjacent to $v$.

A {\it remoteness function} on $K$ for a vertex $v_0$ is a
function $r$ from the set of vertices of $K$ to $\Bbb N \cup
\{0\}$ such that $r^{-1}(0) \subset st(v_0)$.

A function $b$ from the set of adjacency pairs of $K$ to $\Bbb N
\cup\{0\}$ is called a {\it blocking function } for the remoteness
function $r$ if it has the following properties whenever $(X, v)$
is an adjacency pair with $r(v)>0$.
\begin{enumerate}
\item if $b(X, v)=0$, then there exists a vertex $w$ of $K$
which is adjacent to $v$ such that $r(w)<r(v)$ and $(X, w)$ is
also an adjacency pair (see Fig. 3(a)), and
\item if $b(X, v)>0$, then there exist $v'\in X$ and a vertex $w'$ of $K$
which is adjacent to $v'$ such that
 \begin{enumerate}
  \item $r(w')<r(v')$,
  \item every element of $X$ that is adjacent to $v'$
  is also adjacent to $w'$, and
  \item if $Y = (X - \{v'\})\cup\{w'\}$, then $b(Y, v)<b(X, v)$, where
  $Y = (X - \{v'\})\cup\{w'\}$ means remove one instance of $v'$ from $X$
  and add one instance of $w'$ to $X$ (see Fig. 3(b)).
 \end{enumerate}
\end{enumerate}

\begin{figure}
\centering
\includegraphics{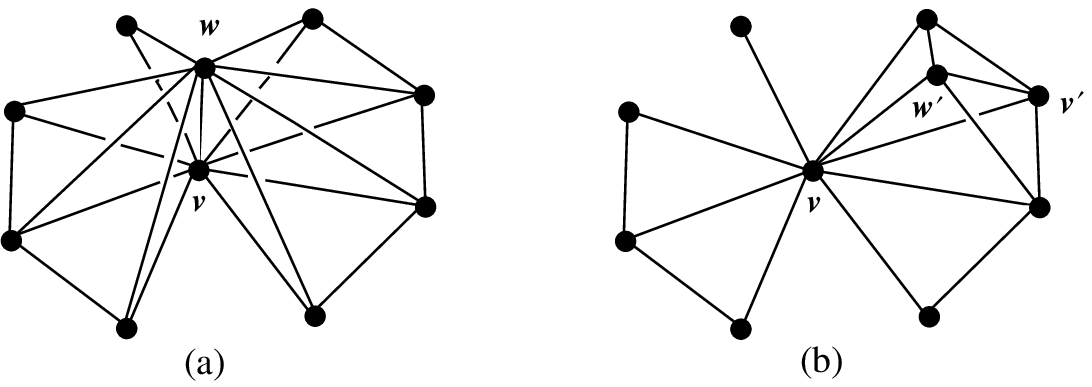}
\caption{}
\end{figure}

A simplicial complex $K$ is said to be {\it flag} if any collection
of $k+1$ pairwise distinct vertices of $K$ spans a $k$-simplex
whenever any two span a 1-simplex. The proof of the following
proposition is based on the proof of Theorem 5.3 in \cite{McC}.

\begin{prop}
Let $K$ be a flag complex with base vertex $v_0$. If $K$ has a
remoteness function $r$ for $v_0$ that admits a blocking function
$b$, then $K$ is contractible.
\end{prop}

\begin{proof}
 It suffices to show that the homotopy groups are all trivial.
 Let $f:\Bbb S^q \rightarrow K$, $q \ge 0$, be a map carrying the base
point of $\Bbb S^q$ to $v_0$. We may assume that $f$ is simplicial
with respect to a triangulation $\Delta$ of $\Bbb S^q$. If
$r(f(u))=0$ for every vertex $u$ of $\Delta$, then the image of
$f$ lies in $st(v_0)$. Since $K$ is a flag complex, $f$ is
null-homotopic.

Let us assume then that there exists a vertex $u$ of $\Delta$ such
that $r(f(u))>0$. Choose such a vertex $u$ so that $r(f(u))=n$ is
maximal among all vertices of $\Delta$. Our goal is to find a
simplicial map $g:\Bbb S^q \rightarrow K$ with respect to some
subdivision $\Delta'$ of $\Delta$ such that:
\begin{itemize}
\item $g$ is homotopic to $f$,
\item $r(g(u))<n$,
\item $r(g(z))=r(f(z))$ for every vertex $z$ of $\Delta$ with $z
\neq u$, and
\item $r(g(u'))<n$ for every vertex $u'$ of $\Delta'-\Delta$.
\end{itemize}
Then, repeating the process for other vertices whose values by $r
\circ f$ are also $n$, we obtain a simplicial map $h$ with respect
to some subdivision of $\Delta$ such that $h$ is homotopic to $f$
and $r(h(u))<n$ for every vertex $u$. We thus complete the proof
inductively.

Let $u_1 , u_2 , \cdots ,u_s $ be the vertices in the link of $u$ in
$\Delta$ and let $f(u)=v$, $f(u_j )=v_j$ for $j=1, \cdots ,s$, and
$X=\{v_1 , v_2 , \cdots ,v_s \}$. Then $X$ is a finite multiset and
$(X, v)$ is an adjacency pair. If $b(X, v)=0$, then there exists a
vertex $w$ in $st(v)$ such that $r(w)<r(v)$ and $(X, w)$ is an
adjacency pair. Define a simplicial map $g:\Bbb S^q \rightarrow K$
with respect to the same triangulation $\Delta$ by $g(u)=w$ and
$g(z)=f(z)$ for all vertices $z \neq u$. Since $K$ is a flag
complex, $g$ is homotopic to $f$ and we have $r(g(u))<n$.

Now suppose $b(X, v)>0$. Then there exist $v_j \in X$, and a vertex
$w_j$ of $K$ which is adjacent to $v_j$ such that (1)
$r(w_j)<r(v_j)$, (2) every element of $X$ that is adjacent to $v_j$
is also adjacent to $w_j$, and (3) $b(Y, v)<b(X, v)$ where $Y = (X -
\{v_j\})\cup\{w_j\}$. Construct a subdivision $\Delta'$ of $\Delta$
by introducing the barycenter $u_j '$ of the simplex $ \langle u,
u_j \rangle $ as a vertex, and replacing each simplex of the form
$\langle u, u_j, z_1, z_2, \cdots, z_r \rangle$ by the two simplices
$\langle u, u_j', z_1, z_2, \cdots, z_r \rangle$ and $\langle u_j',
u_j, z_1, z_2, \cdots, z_r \rangle$. Define a simplicial map
$f':\Bbb S^q \rightarrow K$ with respect to $\Delta'$ by
$f'(u_j')=w_j$ and $f'(z)=f(z)$ for every vertex $z$ of $\Delta$.
Since $K$ is a flag complex, $f'$ is homotopic to $f$. Now $Y$ is
the image of the vertices of the link of $u$ in $\Delta'$, and
$r(f'(u'_j))= r(w_j) < r(v_j)\leq r(v)= n$.

Repeating finitely many times, we obtain a subdivision $\Delta''$ of
$\Delta$ and a simplicial map $f''$ with respect to $\Delta''$ so
that $f''$ is homotopic to $f$ and $b(Y'', v)=0$, where $Y''$ is the
image of the vertices of the link of $u$ in $\Delta''$. Observe that
$r(f''(u''))<n$ for every vertex $u''$ of $\Delta''-\Delta$ and
$f''(z)=f(z)$ for every vertex $z$ of $\Delta$. Since $b(Y'', v)=0$,
we obtain a simplicial map $g:\Bbb S^q \rightarrow K$ with respect
to $\Delta''$ as above such that $g$ is homotopic to $f$,
$r(g(u))<n$, and $r(g(u''))<n$ for every vertex $u''$ of
$\Delta''-\Delta$.
\end{proof}

\section{The complex of primitive disks}

The {\it disk complex} $D(V_g)$ of a handlebody $V_g$ of genus
$g$, for $g\geq 2$, is a simplicial complex defined as follows.
The vertices of $D(V_g)$ are isotopy classes of essential disks in
$V_g$, and a collection of $k+1$ vertices spans a $k$-simplex if
and only if it admits a collection of representative disks which
are pairwise disjoint. When $V_g$ is a handlebody in a Heegaard
splitting $(V_g , W_g ; \Sigma_g)$ of $\Bbb S^3$, the {\it complex
of primitive disks} $P(V_g)$ of $V_g$ is defined to be the full
subcomplex of $D(V_g)$ spanned by vertices whose representatives
are primitive disks in $V_g$. As before, we write $V$ and $(V , W
; \Sigma)$ for $V_2$ and $(V_2 , W_2 ; \Sigma_2)$ respectively.
Notice that $D(V)$ and $P(V)$ are 2-dimensional.

It is a standard fact that any collection of isotopy classes of
essential disks in $V_g$ can be realized by a collection of
representative disks that have pairwise minimal intersection. One
way to see this is to choose the disks in their isotopy classes so
that their boundaries are geodesics with respect to some hyperbolic
structure on the boundary surface $\Sigma_g$ and then remove simple
closed curve intersections of the disks by isotopy. In particular,
for a collection $\{v_0 , v_1 , \cdots , v_k \}$ of vertices of
$D(V_g)$, if $v_i$ and $v_j$ bound a 1-simplex for each $i<j$, then
$\{v_0 , v_1 , \cdots , v_k \}$ is realized by a collection of
pairwise disjoint representatives. Thus we have

\begin{lemma}
$D(V_g)$ is a flag complex, and consequently any full subcomplex
of $D(V_g)$ is also a flag complex.
\end{lemma}

\begin{theorem}
If $K$ is a full subcomplex of $D(V_g)$ satisfying the following
condition, then $K$ is contractible.

\begin{itemize}
 \item Suppose $E$ and $D$ are any two disks in $V_g$ which represent vertices of $K$ such that $E
\cap D \neq \emptyset$. If $C\subset D$ is a disk cut off from $D$
by an outermost arc $\alpha$ of $D \cap E$ in $D$ such that $C\cap
E= \alpha$, then at least one of the disks obtained from surgery
on $E$ along $C$ also represents a vertex of ~ $K$.
\end{itemize}
\end{theorem}

\begin{proof}
Since $K$ is a flag complex, by Lemma 4.1, it suffices to find a
remoteness function that admits a blocking function as in
Proposition 3.1. Fix a base vertex $v_0$ of $K$. Define a
remoteness function $r$ on the set of vertices of $K$ by putting
$r(w)$ equal to the minimal number of intersection arcs of disks
representing $v_0$ and $w$.

Let $(X, v)$ be an adjacency pair in $K$ where $r(v)>0$ and
$X=\{v_1, v_2, \cdots , v_n\}$. Choose representative disks  $E,
E_1, \cdots ,E_n $ and $D$ of $v, v_1, \cdots , v_n $ and $v_0$
respectively so that they have transversal and pairwise minimal
intersection. Since $r(v)>0$, so $D\cap E \neq \emptyset$. Let
$C\subset D$ be a disk cut off from $D$ by an outermost arc
$\alpha$ of $D\cap E$ in $D$ such that $C\cap E= \alpha$. Observe
that each $E_i$ is disjoint from the arc $\alpha$.

Let $ b_0 = b_0 (E, E_1, \cdots , E_n, D)$ be the minimal number
of arcs $\{C \cap E_1 \} \cup \{C \cap E_2 \} \cup \cdots \cup \{C
\cap E_n \}$ in $C$ as we vary over such disks $C$ cut off from
$D$. Define $b=b(X, v)$ to be the minimal number $b_0$ as we vary
over such representative disks of $v, v_1, \cdots , v_n $ and
$v_0$. We verify that $b$ is a blocking function for the
remoteness function $r$ as follows.

Suppose, first, $b(X, v)=0$. There exist then representative disks
$E', E'_1, \cdots ,E'_n $ and $D'$ of $v, v_1, \cdots , v_n $ and
$v_0$ respectively and exists, as above, a disk $C'\subset D'$ cut
off from $D'$ so that $C' \cap (E'_1 \cup \cdots \cup
E'_n)=\emptyset$. By the assumption, a disk obtained from surgery
on $E'$ along $C'$ represents a vertex $w'$ of $K$ again, and $(X,
w')$ is an adjacent pair. We have $r(w')<r(v)$ since the arc $C'
\cap E'$ no longer counts.

Next, suppose $b(X, v)>0$. Choose then representative disks $E'',
E''_1, \cdots ,E''_n $ and $D''$ of $v, v_1, \cdots , v_n $ and
$v_0$ respectively, and, similarly, a disk $C''\subset D''$ cut
off from $D''$ so that they realize $b(X, v)$. Choose an outermost
arc $\beta$ of $C'' \cap E''_k$ in $C''$ that cuts off a disk
$F''\subset C''$ such that (1) $F'' \cap E''_k = \beta$, (2) $F''$
is disjoint from the arc $C'' \cap E''$, and (3) $F''$ contains no
arc of intersection of any $E''_i$ which is disjoint from $E''_k$
(note that $\beta$ may still intersect an arc of some $C'' \cap
E''_i$). Then a disk obtained from surgery on $E''_k$ along $F''$
represents a vertex $w''$ of $K$ by the assumption. By the
construction, every element of $X$ that is adjacent to $v_k$ is
also adjacent to $w''$. We have $r(w'')<r(v_k)$ and $b(Y, v)<b(X,
v)$, where $Y = (X - \{v_k\})\cup\{w''\}$, since $\beta$ no longer
counts.
\end{proof}

Theorem 4.2 shows that $D(V_g)$ is contractible since a disk
obtained from surgery on an essential disk is also essential. Now we
return to the genus two case.

\begin{cor}
$P(V)$ is contractible, and the link of any vertex of $P(V)$ is
contractible.
\end{cor}

\begin{proof}
From Theorem 2.3, $P(V)$ and the link of a vertex $v$ in $P(V)$
satisfy the condition of Theorem 4.2. For the case of the link of
$v$, observe that one of the primitive disks from surgery
represents $v$ again, but the other one represents a vertex of the
link.
\end{proof}

Let $P(V)'$ be the first barycentric subdivision of $P(V)$ and let
$T$ be the graph formed by removing all of the open stars of the
vertices of $P(V)$ in $P(V)'$. Each edge of $T$ has one endpoint
which represents the isotopy class of a reducing triple (the
barycenter of a 2-simplex of $P(V)$) and one endpoint which
represents the isotopy class of a reducing pair (the barycenter of
an edge of that 2-simplex) in that reducing triple. By Corollary
4.3, we see that $T$ is contractible. Thus we have

\begin{cor}
The graph $T$ is a tree.
\end{cor}

\section{A finite presentation of ${\mathcal H_2}$}

In this section, we will use our description of $T$ to recover
Akbas's presentation of ${\mathcal H_2}$ given in \cite{Ak}. The
tree $T$ is invariant under the action of $\mathcal H_2$. In
particular, $\mathcal H_2$ acts transitively on the set of
vertices of $T$ which are represented by reducing triples. The
disks in a reducing triple are permuted by representative
homeomorphisms of $\mathcal H_2$. It follows that the quotient of
$T$ by the action of $\mathcal H_2$ is a single edge. For
convenience, we will not distinguish disks and homeomorphisms from
their isotopy classes in their notations.

\begin{figure}
\centering
\includegraphics{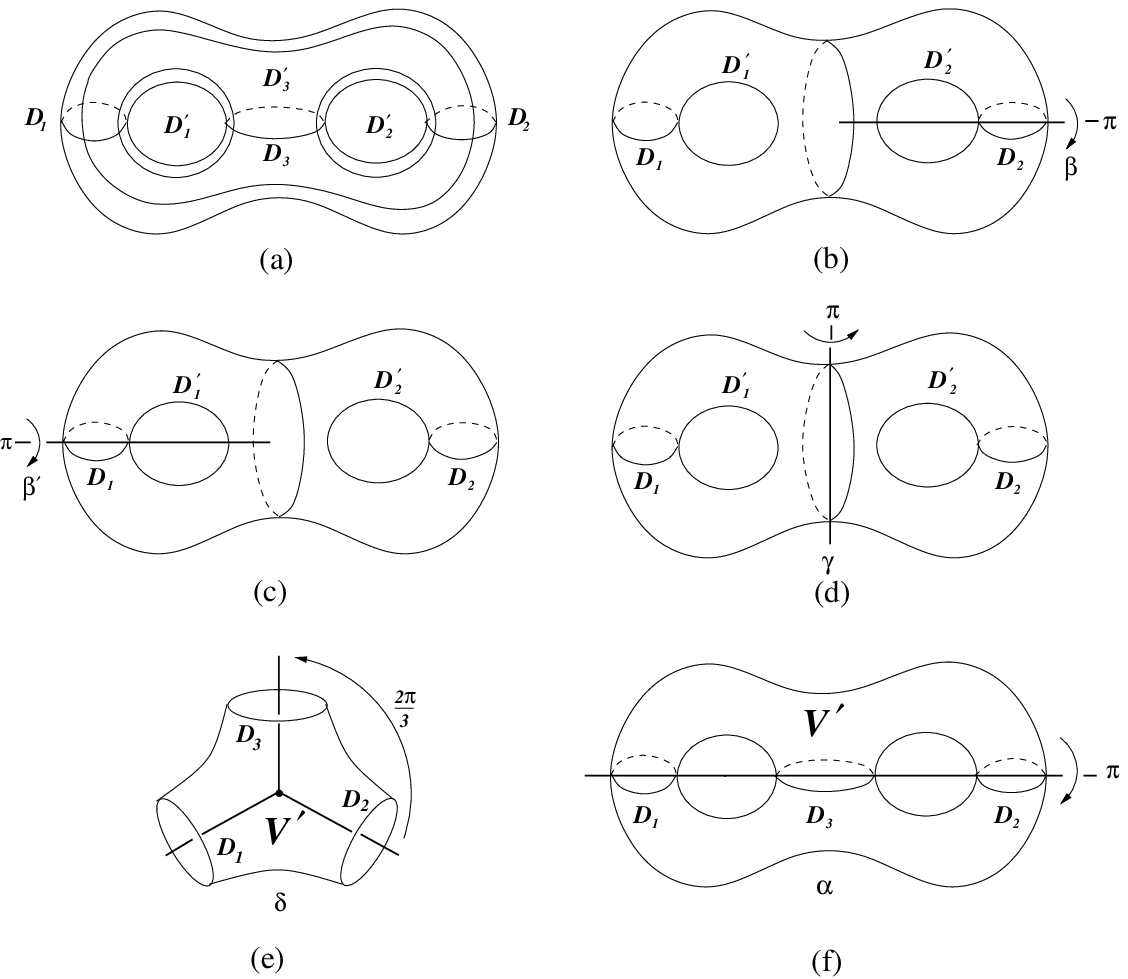}
\caption{}
\end{figure}

Fix two vertices $P=\{D_1, D_2\}$ and $Q=\{D_1, D_2, D_3\}$ that are
endpoints of an edge $e$ of $T$. Denote by $\mathcal H _P $,
$\mathcal H _Q$ and $\mathcal H_ e$ the subgroups of $\mathcal H_2$
that preserve $P$, $Q$ and $e$ respectively. By the theory of groups
acting on trees due to Bass and Serre \cite{S}, $\mathcal H_2$ can
be expressed as the free product of $\mathcal H _P $ and $\mathcal H
_Q$ with amalgamated subgroup ~$\mathcal H_ e$.

It is not difficult to describe these subgroups of $\mathcal H_2$.
We sketch the argument as follows. Let $D_1'$ and $D_2'$ be the
disjoint dual disks of $D_1$ and $D_2$ respectively, and let $D_3'$
be the dual disk of both $D_1$ and $D_2$ which is disjoint from
$D_1'$, $D_2'$ and $D_3$ (see Fig. 4(a)). The reducing triple
$\{D_1', D_2', D_3'\}$ of $W$ is uniquely determined by Lemma 2.2.
Denote $\partial D_i$ and $\partial D_i'$ by $d_i$ and $d_i'$, for
$i\in\{1, 2, 3\}$, respectively.

Consider the group $\mathcal H _P $. The elements of $\mathcal H _P
$ also preserve $D_1' \cup D_2'$. Let $\mathcal H_P'$ be the
subgroup of elements of $\mathcal H_P$ which preserve each of $D_1$
and $D_2$. Since the union of $\Sigma$ with the four disks $ D_i$
and $ D_i'$, for $i\in\{1, 2\}$, separates $S^3$ into (two) 3-balls,
$\mathcal H_P'$ can be identified with the group of isotopy classes
of orientation preserving homeomorphisms of the annulus obtained by
cutting $\Sigma$ along $d_1 \cup d_1' \cup d_2 \cup d_2'$, which
preserve each of $d_1$ and $d_2$. This group is generated by two
elements $\beta$ and $\beta'$ ($\pi$-twists of each boundary circle)
and has relations $(\beta \beta')^2 =1$ and
$\beta\beta'=\beta'\beta$. As generators of $\mathcal H _P' $,
$\beta$ and $\beta'$ are shown in Fig. 4(b) and Fig. 4(c).

Let $\gamma$ be the order two element of $\mathcal H_P$ which
interchanges $D_1$ and $D_2$, as shown in Fig. 4(d). Then
$\mathcal H_P$ is an extension of $\mathcal H_P'$ by $\langle
\gamma \rangle$ with relations $\gamma \beta \gamma =
{\beta'}^{-1}$ and $\gamma \beta' \gamma = {\beta}^{-1}$.

Next, let $V'$ be a 3-ball cut off from $V$ by $D_1 \cup D_2 \cup
D_3$, and let $\mathcal H_Q'$ be the subgroup of elements of
$\mathcal H_Q$ which preserve $V'$. Since the union of $\Sigma$ with
the six disks $ D_i$ and $ D_i'$, for $i\in\{1, 2, 3\}$, separates
$S^3$ into (four) 3-balls, $\mathcal H_Q'$ can be identified with
the group of isotopy classes of orientation preserving
homeomorphisms of the 3-holed sphere $\partial V' \cap \partial V$,
which preserve $d_1 \cup d_2 \cup d_3$ and $d_1' \cup d_2' \cup
d_3'$ respectively. This group is a dihedral group generated by
$\gamma $ and $\delta$, where $\delta$ is the order three element
shown in Fig. 4(e).

Let $\alpha$ be the order two element of $\mathcal H_Q$ which
preserves each $D_1$, $D_2$ and $D_3$ but takes $V'$ to another
component cut off from $V$ by $D_1 \cup D_2 \cup D_3$, as shown in
Fig. 4(f). Then $\mathcal H_Q$ is an extension of $\mathcal H_Q'$ by
$\langle \alpha \rangle$ with relations $\alpha \gamma \alpha =
\gamma$ and $\alpha \delta \alpha = \delta$. Finally, it is easy to
see that $\mathcal H_ e$ is generated by $\gamma$ and $\alpha$.
Observing that $\beta'=\gamma {\beta}^{-1} \gamma$ and
$\alpha=\beta\beta'$, we get a finite presentation of $\mathcal H_2$
with generators $\beta$, $\gamma$, and $\delta$.

\section{The Scharlemann-Akbas Tree}

Scharlemann \cite{Sc} constructed a connected simplicial 2-complex
$\Gamma$ which deformation retracts to a certain graph
$\tilde{\Gamma}$ on which $\mathcal H_2$ acts, and gave a finite
generating set of $\mathcal H_2$. In this section, we describe
$\Gamma$ and $\tilde{\Gamma}$ briefly and compare them with $P(V)$
and $T$. For details about $\Gamma$ and $\tilde{\Gamma}$, we refer
the reader to \cite{Sc} and \cite{Ak}.

Let $P$ and $Q$ be reducing spheres for a genus two Heegaard
splitting $(V, W; \Sigma)$ of $\Bbb S^3$. Then $P\cap V$ and
$Q\cap V$ are essential separating disks in $V$ and each of them
cuts $V$ into two unknotted solid tori. Define the {\it
intersection number} $P \cdot Q$ to be the minimal number of arcs
of $P \cap Q \cap V$ up to isotopy. Then we observe that $P \cdot
Q \ge 2$ if $P\cap V$ and $Q\cap V$ are not isotopic to each other
in $V$.

The simplicial complex $\Gamma$ for $(V, W; \Sigma)$ is defined as
follows. The vertices are the isotopy classes of reducing spheres
relative to $V$ and a collection $P_0, P_1, \cdots , P_k$ of $k+1$
vertices bounds a $k$-simplex if and only if $P_i \cdot P_j =2$ for
all $i<j$. It turns out that $\Gamma$ is a 2-complex and each edge
of $\Gamma$ lies on a single 2-simplex. Thus $\Gamma$ deformation
retracts naturally to a graph $\tilde{\Gamma}$ in which each
2-simplex in $\Gamma$ is replaced by the cone on its three vertices,
and Akbas \cite{Ak} showed that $\tilde{\Gamma}$ is a tree.

\begin{figure}
\centering
\includegraphics{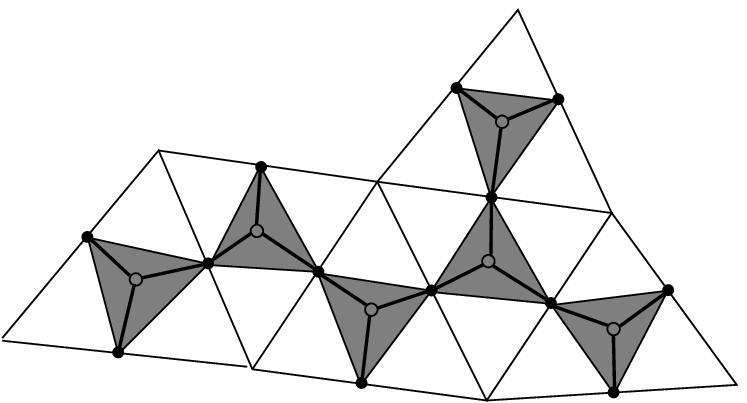}
\caption{}
\end{figure}

We observe that the reducing spheres for $(V, W;\Sigma)$
correspond exactly to the reducing pairs in $V$ up to isotopy.
That is, a reducing sphere cuts $V$ into two unknotted solid tori,
whose meridian disks are primitive in $V$. Conversely, a reducing
pair has a unique  pair of disjoint dual disks, as in Lemma 2.2,
and the corresponding reducing sphere is the boundary of a small
regular neighborhood of the union of either of the disks with its
dual.

It is routine to check that two reducing spheres represent two
vertices connected by an edge in $\Gamma$ if and only if the
corresponding reducing pairs are contained in a reducing triple.
So the correspondence from reducing spheres to reducing pairs
determines a natural embedding of $\Gamma$ into $P(V)$ so that the
image of $\tilde{\Gamma}$ is identified with the tree $T$ (see
Fig. 5). Therefore, Corollary 4.4 can be considered as an
alternate proof that $\tilde{\Gamma}$ is a tree.

\bibliographystyle{amsplain}

\end{document}